\documentclass[12pt]{article}

\usepackage{amssymb,amsmath,amsthm,bbm,epsfig}
\usepackage{graphicx}
\usepackage{psfrag}

\theoremstyle{definition}
\newtheorem{thm}{Theorem}[section]
\newtheorem{lem}{Lemma}[section]
\newtheorem{cor}{Corollary}[section]

\newtheorem*{remark}{Remark}

\newcommand{\C}{\mathbbm{C}}   
\newcommand{\R}{\mathbbm{R}}                 %%%% R
\newcommand{\N}{\mathbbm{N}} 
\newcommand{\Z}{\mathbbm{Z}}                 %%%% Z

\newcommand{\D}{\mathbb{D}} 
\newcommand{\opdisk}{\mathbb{U}} 

\newcommand{\p}{\mathbbm{P}}
\newcommand{\pr}[1]{\mathbb{P}\left\{#1\right\}}
\newcommand{\pro}[2]{\mathbb{P}^{#1}\left\{#2\right\}}

\newcommand{\ex}[1]{\mathbb{E}\left[#1\right]}

\newcommand{\im}[1]{\text{Im}(#1)}
\newcommand{\re}[1]{\text{Re}(#1)}

\newcommand{\gam}{\gamma}
\newcommand{\eps}{\epsilon}

\newcommand{\bigo}[1]{\mathcal{O}\left(#1\right)}

\renewcommand{\S}[3]{\displaystyle{\sum_{#1}^{#2}{#3}}} 

\newcommand{\bd}[1]{\partial{#1}}

\renewcommand{\dh}{\mathcal{H}}

\newcommand{\h}{\mathbbm{H}}

\begin{document}

%\begin{center}
%\underline{\huge Random Walk Holes via Imbedding}
%\end{center}

\title{Some Estimates for Planar Random Walk and Brownian Motion}
\author{Christian Bene\v{s}\\ Tufts University}
\maketitle

\begin{abstract}
 The purpose of this note is to collect in one place a few results about simple
 random walk and Brownian
 motion which are often useful. These include standard results such as Beurling estimates, large deviation estimates,
 and a method for coupling the two processes, as well as solutions to the discrete
 Dirichlet problem in various domains which, to the author's knowledge, have
 not been published anywhere. The main focus is on the two-dimensional
 processes.
\end{abstract}

\section{Introduction and Definitions}\label{IntroDefs}

In the study of simple random walk or standard Brownian motion, some estimates are
striking by their ubiquity. The following pages are a collection of
results which are frequently needed when dealing with these processes, and of some more specific
estimates.

$B$ will denote planar standard Brownian motion started at the origin and $S$ will be planar simple random walk, 
defined by $S(0)=(0,0)$ and for $n\in\N$, by $S(n) = \sum_{k=1}^{n}X_k$,
where $\{X_k\}_{k\in\N}$ are independent random vectors satisfying
$\pr{X_k = \pm e_i} = \frac{1}{4}, i=1,2$,
where $e_1=\langle 1,0\rangle$ and $e_2=\langle 0,1\rangle$. 
We will also think of planar simple random walk $S$ as being a
continuous process, that is, for positive real times $t$, we let $S(t)$ be
the linear interpolation of the walk's position at the surrounding
integer times: For all real $t\geq 0$, 
\begin{equation}\label{03-28-05}
S(t) = S([t])+(t-[t])(S([t]+1)-S([t])),
\end{equation}
where $[t]$ denotes the integer part of $t$. For any real numbers
$0\leq a\leq b$, we will write
$S[a,b]:=\{S(t)\}_{a\leq t\leq b}$, and use the same notation for $B$. 
If we consider $B$ or $S$ in dimensions other than 2, we will say so explicitly. The measure associated with either of the processes, started at $x$ will be denoted $\p^x$ and we will write $\p$ for $\p^0$. It will be clear from context which is the concerned process.

In what follows, all multiplicative constants will be denoted by
$K, K_1$, or $K_2$. It will be understood that they may be
different from one line to the next. The letters $r, s, t$ will be used to denote real numbers, while
$i, j, k, l, m, n$ will be integers. 

%Points of the
%complex plane $\C$ will be represented by the letters $u, v, w, z$.

The symbols ${\cal O}$, $o$, $\sim$, and $\asymp$ will mean the following: for two
functions $f$ and $g$, $f(x) = \bigo{g(x)}$ if there exists a constant $K$ such that
$f(x)\leq Kg(x)$ for all $x$ large enough, 
$f(x)=o(g(x))$ if
$\lim_{x\to\infty} f(x)/g(x)=0$, $f(x)\sim
g(x)$ if $\displaystyle{\lim_{x\to\infty} f(x)/g(x)=1},$ and $f(x)
\asymp g(x)$ if there exists a constant $K>0$ such that for all $x$ large enough
$\frac{1}{K}g(x)\leq f(x) \leq Kg(x)$. In some cases, the same notation will be used to describe limiting behavior close to 0 and the context will leave no doubt as to what is meant.
%We will often be interested in the
%dependence of $K$ on variables involved in our equations and will
%refer to it as the constant of the $\mathcal{O}$. 
%We say that a
%function $f(x)$ {\it decays rapidly} if for every $r\geq 0, f(x) = o(x^{-r})$.
% The boundary
%of a set $A\subset\C$ will be denoted by $\bd{A}$, its area by $|A|$,
%and its diameter $\text{diam}(A)=\sup_{x,y\in A}|x-y|$. For any sets
%$A\subset F$, $A^c = F\setminus A$ will be the complement of $A$ in
%$F$. It will always be clear from context what $F$ is meant to be. 
The Euclidean norm of a point $x$ in $\R$ or $\C$ is $|x|$.  The distance between two sets $A,B\subset\C$ is $d(A,B)=\inf_{x\in A, y\in B} |x-y|$ and the diameter of $A$ is $\text{diam}(A)=\sup_{w,z\in A}|w-z|$. The boundary of a discrete set $D\subset \Z^2$ is defined to be $\{w\in \Z^2:w\not\in D; \exists z\in D \text{ with } |z-w|=1\}$.

%Unless stated otherwise, $B=(B(t))_{t\geq 0}$ will denote standard planar Brownian motion and
%$S=(S(t))_{t\geq 0}$ will stand for planar simple random walk. 

\section{From One to Two Dimensions}\label{1-2}

Planar Brownian motion $B$ can be defined as a couple
$(B^1(t),B^2(t))$ of independent realizations of one-dimensional
Brownian motion. This representation is often convenient when trying to obtain
estimates for the planar process from estimates for the linear process.
One does not usually think of planar random walk in
the same way, although an analogous definition is possible. Since this will be useful in the next section, we quickly give
this simple construction here:

The idea is to note that if we take two independent random walks $S^1$
and $S^2$ on appropriately scaled and rotated versions of $\Z$, then $(S^1,S^2)$ is a
simple random walk in $\Z^2$.
\begin{figure}
         \centering
         \includegraphics[width=.5\textwidth]{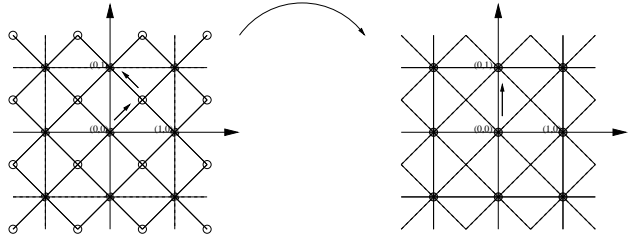}
         \caption{Two independent one-dimensional simple random walks on diagonal lattices give a simple random walk in $\Z^2$ (black dots).}
         \label{Fig3}
\end{figure}
More precisely, for $k\geq 1, j=1,2$, let $X_k^j$ be independent random vectors with distribution
$$\pr{X_k^1=\pm \frac{1}{\sqrt{2}}e^{i\pi/4}} = \pr{X_k^2=\pm\frac{1}{\sqrt{2}}e^{i3\pi/4}} = \frac{1}{2}.$$
Then if we let $S^j(n)=\S{k=1}{n}{X_k^j}, S^1(n)$ and $S^2(n)$  are
independent simple random walks on
$\frac{1}{\sqrt{2}}e^{i\pi/4}\cdot\Z =
\{\frac{l}{\sqrt{2}}e^{i\pi/4}:l\in\Z\}$ and
$\frac{1}{\sqrt{2}}e^{i3\pi/4}\cdot\Z$, respectively, and $S(n) =
(S^1(n),S^2(n))$ is a simple random walk in $\Z^2$. See Figure \ref{Fig3}.

\section{Large deviations}\label{LargeDev}

In time $n$, planar Brownian motion $B$ and planar simple random walk
$S$
are expected to reach a distance of order $\sqrt{n}$, as shown in the first lemma below. In this
section we give bounds for the likelihood that they behave
unusually (that is, reach a distance which is much greater than
$\sqrt{n}$ or remain in a disk which has a
radius much smaller than $\sqrt{n}$). We only give upper bounds for
these probabilities, except in the case of Brownian motion traveling much beyond distance $\sqrt{n}$, where deriving a lower bound requires little
work. The corresponding results for random walk are remarkably more
difficult to obtain. 

%For $x \in \C$ and $d \in \R_+$, let 
%$$D(x,d):=\{y \in \C:|y-x|\leq d\}, \;\;\;\;\; \tilde{D}(x,d)=\{z \in \Z^2:|z-x|\leq d\}.$$
%$$\tilde{\xi}_r = \inf\{k \geq 0:|S_k - S_0| \geq r\}, \hspace{3pc} \xi_r = \inf\{t \geq 0:|B(t) - B_0| \geq r\}. $$
%When we write ``$\forall r\geq c, X(r) \asymp Y(r)$'', we mean: $\exists$ constants $C_1, C_2$ such that $\forall r\geq c, C_1Y(r) \leq X(r) \leq C_2Y(r)$.

%We first show how to find for a planar symmetric process $X$ (continuous
%or discrete) $\pr{sup_{0\leq t\leq n}|X(t)|\geq a}$ from $\pr{X(n)\geq
%  a}$.

%\begin{def} We say that a discrete planar Markov process $X$ is {\bf
%    symmetric} if for every $\vec{a}\in\R^2$ its transition
%    probabilities satisfy $\pr{S(n+1)-S(n) = \vec{a}} =
%    \pr{S(n+1)-S(n) = -\vec{a}}$. A continuous planar Markov process
%    is symmetric if 

\begin{lem}\label{meansquaredistance}
\begin{eqnarray*}
& (a) & \ex{|B(n)|^2}=2n.\\
& (b) & \ex{|S(n)|^2}=n.
\end{eqnarray*}
\end{lem}

\begin{proof}
(a) The transition density for $B$ is $p_t(x,y) = \frac{1}{2\pi t}e^{-|y-x|^2/2t}$, so integration in polar coordinates yields
$$\ex{|B(n)|^2} = \int_0^{2\pi}\int_0^\infty \frac{1}{2\pi n}e^{-r^2/2n}r^3\,dr\,d\theta = 2n.$$ 

(b) It is easy to check that $M_n=|S_n|^2-n$ is a martingale, so $\ex{M_n}=\ex{M_0}=0$.
\end{proof}

\begin{lem}\label{reflection}
 If $X$ is planar random walk or Brownian motion, then for any $a\in\R_+, n\in\N$,
$$\pr{|X(n)|\geq a}\leq \pr{\sup_{0\leq s\leq n}|X(s)|\geq a} \leq 2\pr{|X(n)|\geq
  a}.$$
\end{lem}

\begin{proof}
The first inequality is obvious. The second is a consequence of  
\begin{eqnarray*}
\pr{\sup_{0\leq s\leq n}|X(s)|\geq a} & \leq & \pr{\sup_{0\leq
    s\leq n}|X(s)|\geq a ; |X(n)|\geq a}\\
& & \hspace{7.15pc} + \pr{\sup_{0\leq
    s\leq n}|X(s)|\geq a ; |X(n)| < a}\\
& \leq & \pr{|X(n)|\geq a} + \frac{1}{2}\pr{\sup_{0\leq s\leq n}|X(s)|\geq a},
 \end{eqnarray*}
where the last inequality follows from the strong Markov property.
\end{proof}

\subsection{Brownian motion}\label{bm}

Lemma \ref{bmtoofar} gives a two-sided estimate for the probability that Brownian
motion goes farther than expected, while Lemma \ref{bmtooclose} provides an upper bound for the probability that it remains in a ball that is much smaller than expected.

\begin{lem}\label{bmtoofar} If $B$ is planar standard Brownian motion,
\begin{eqnarray*}
& (a) & \pr{|B(n)|\geq r\sqrt{n}} = e^{-r^2/2}.\\
& (b) & \pr{\sup_{0\leq t \leq n}|B(t)| \geq r\sqrt{n}} \asymp  e^{-r^2/2}.
\end{eqnarray*}
\end{lem}

\begin{proof}
$(a)$ This is a straightforward computation using the transition density
 for $B$. $(b)$ follows
 from $(a)$ and Lemma \ref{reflection}.

%If we write $B(t) = (B^1(t),B^2(t))$, Brownian scaling, and the reflection principle give\\

%$\begin{array}{lll}
%\displaystyle{ \quad \pr{B_1^1 \geq r} = \mathbb{P}\{B_n^1 \geq rn^{1/2}\}}&  & \\
%  & \leq \displaystyle{\pr{\sup_{0\leq t \leq n}|B(t)|\geq rn^{1/2}}} & \\
% & & \leq \displaystyle{2\pr{\sup_{0\leq t \leq n}|B^1(t)|\geq r\sqrt{n/2}}}\\
% & & \leq \displaystyle{4\pr{\sup_{0\leq t \leq n}B^1(t)\geq r\sqrt{n/2}}}\\
% & & = \displaystyle{2\pr{B_n^1\geq r\sqrt{n/2}}}\\
% & & = \displaystyle{2\pr{B_1^1\geq r/\sqrt{2}}}.
%\end{array}$\\
%If $\phi$ is the normal density and $r>0$, then
%$$\frac{1}{r + r^{-1}}\phi(r) \leq \pr{B_1^1 \geq r} \leq \frac{1}{r}\phi(r),$$
%from which the result follows immediately
\end{proof}

\begin{remark} The following is an equivalent formulation of part $(b)$
 of the lemma: If $T_R = \inf\{t\geq 0:|B(t)|\geq R\}$,
$$ \pr{T_{r\sqrt{n}} \leq n}\asymp e^{-r^2/2}.$$
%\asymp r^{-1}\exp\{-\frac{r^2}{2}\}.$$
\end{remark}

\begin{lem}\label{bmtooclose} There exists a constant $K > 0$ such that for all $n\in\N, r\geq 1$,
$$\pr{\sup_{0\leq t \leq n}|B(t)| \leq r^{-1}n^{1/2}} \leq \exp\{-K r^2\},$$
where $B$ is a planar Brownian motion.
\end{lem}

\begin{proof}
For $l,n \in\N, r\geq 1$, we define $$I_l = I_l(r,n) = [(l-1)\frac{n}{r^2},l\cdot
  \frac{n}{r^2}].$$ Then, if $k = [r^2]$, where $[\cdot]$ denotes the
  integer part, $$\displaystyle{\cup_{l=1}^k I_l \subset [0,n]}.$$ A simple geometric argument, the Markov property, and Brownian scaling give:
\begin{eqnarray*}
\pr{\sup_{0\leq t \leq n}|B(t)| \leq r^{-1}n^{1/2}} & \leq & \prod_{l=1}^{k}\pr{\sup_{t\in I_l}|B(t) - B\left((l-1)\frac{n}{r^2}\right)| \leq 2r^{-1}n^{1/2}}\\
& \leq & \prod_{l=1}^{k}\pr{\sup_{0\leq t \leq 1}|B(t)|\leq 2}\\
& \leq & \rho^k \leq \exp\{-K r^2\},
\end{eqnarray*}
where $\rho < 1$ is independent of $r$ and $n$, and $K=-\frac{\ln \rho}{2} > 0$.

\end{proof}

\subsection{Random walk}\label{rw}

As mentioned above, deriving the same type of estimates for random walk
is more involved, since we do not have as nice a transition density to work with. We give in Theorem \ref{LCLTthm} a version of the local central limit theorem which is sharper than, say, in \cite{greenbook} for points far away from the origin and allows us to obtain a sharp large deviations estimate in the one-dimensional case. We then use that estimate to find a bound in the planar case. 

We first prove under slightly
 more general assumptions than for simple random walk a weaker result for random walk which will be useful in Section \ref{Skorokhod}.

\begin{lem}\label{expboundrv} Suppose $\{X_i\}_{i\geq 1}$ are independent, identically distributed random variables with mean 0 such that for some $a>0$, the moment generating function
\begin{equation}\label{exprv} 
M(t)=\ex{e^{tX_1}}<\infty \text{ for }|t|\leq a.
\end{equation}
Then there exists a constant $K$ depending on $a$ and the distribution of $X_1$ such that if $S(n)=\S{i=1}{n}{X_i}$, for all $n\geq 1, r>0$, 
$$\pr{|S(n)|\geq r\sqrt{n}}\leq Ke^{-ar}$$
and 
$$\pr{\max_{1\leq k \leq n}|S(k)|\geq r\sqrt{n}}\leq 2Ke^{-ar}$$
\end{lem}

\begin{proof} Choose a distribution $X_1$ and an $a$ for which (\ref{exprv}) holds. It suffices to show that there is a $K_1>0$ such that $\pr{S(n)\geq r\sqrt{n}}\leq K_1e^{-ar}$. By Markov's inequality,
$$\pr{\frac{S(n)}{\sqrt{n}}\geq r}\leq \frac{\ex{\exp\{aS(n)/\sqrt{n}\}}}{e^{ar}}.$$
But by expanding $\ex{e^{tX_1}}$ about 0, we get for $|t|\leq a$,
$$M(t) = 1+\frac{\ex{X_1^2}}{2}t^2+\bigo{t^3},$$
so that we can find a constant $K_2$ such that for all $n\geq 1$,
$$M(\frac{a}{\sqrt{n}})\leq 1+\frac{K_2}{n}.$$
Thus, 
$$\ex{\exp\{aS(n)/\sqrt{n}\}} = M\left(\frac{a}{\sqrt{n}}\right)^n \leq \left(1+\frac{K_2}{n}\right)^n \leq K_1,$$ 
where $K_1 = e^{K_2}$. This implies that 
$$\pr{\frac{S(n)}{\sqrt{n}}\geq r} \leq K_1e^{-ar}.$$
The other part follows from Lemma \ref{reflection}.
\end{proof}

%The next Corollary is an immediate consequence. 

%\begin{cor}\label{rwtoofar} Under the hypotheses of Lemma \ref{expboundrv}, there exists a constant $K>0$ such that for all $n\geq 1, r>0$,
%$$\pr{\max_{1\leq k \leq n} |k|\geq r\sqrt{n}}\leq Kne^{-ar}.$$
%\end{cor}

The following particular case will be of special interest when we
consider Skorokhod embedding in Section \ref{Skorokhod}:

\begin{cor}\label{bmstoppingtimes} For one-dimensional standard Brownian motion started at 0, define $T_1=\inf\{t\geq 0:|B(t)|=1\}$ and for $j\geq 2, T_j=\inf\{t\geq T_{j-1}:|B(t)-B(T_{j-1})|=1\}.$ Then there exist constants $K, K_1 >0$ such that for all $n\geq 1, r>0$,
$$\pr{\max_{1\leq k\leq n}|T_k - k|\geq r\sqrt{n}}\leq Kne^{-K_1r}.$$
\end{cor}
\begin{proof} The fact that $\ex{T_1}=1$ is well known and since for all $j\geq 1, T_j - T_{j-1} \stackrel{\mathcal{D}}{=} T_1$, it suffices, by Lemma \ref{expboundrv} to show that there is an $a>0$ such that $\ex{e^{a(T_1-1)}}<\infty$. It follows from
$$\pr{T_1\geq k+1 | T_1\geq k}\leq \pr{|B(k+1)-B(k)|\leq 2} = \rho < 1$$
that $\pr{T_1\geq k}\leq \rho^k$, and it suffices to choose $a<\ln (\rho^{-1})$ to ensure that $\ex{e^{aT_1}}<\infty$.
\end{proof}

To find a better bound for $\pr{|S(n)|\geq r\sqrt{n}}$ when $S(n)$ is
planar simple random walk, we will derive a bound in the one-dimensional
case via a precise version of the Local Central Limit Theorem and
use it to find a bound in the two-dimensional case.

\begin{thm}[Local Central Limit Theorem]\label{LCLTthm}
If $S$ is one-dimensional simple random walk, then for every $n\geq 1,
|k|\leq n$, 
$$\pr{S(2n)=2k} = \sqrt{\frac{1}{\pi n}}\exp(-\phi(n,k))\left(1+\bigo{\frac{1}{n}}+\bigo{\frac{k^2}{n^2}}\right),$$
where $\phi(n,k) =
    \S{l=1}{\infty}{\frac{1}{l(2l-1)}\frac{k^{2l}}{n^{2l-1}}}$.\\ 
In particular, for every $N\geq 2, a<\frac{2N-1}{2N}$, and $|k|\leq n^a$, 
\begin{equation}\label{LCLTeq1}
\pr{S(2n)=2k} \sim \sqrt{\frac{1}{\pi
    n}}\exp\left(-\S{l=1}{N-1}{\frac{1}{l(2l-1)}\frac{k^{2l}}{n^{2l-1}}}\right).
\end{equation}
Moreover, there exists a constant $K>0$ such that for all $n\in\N, k\in\Z$, 
\begin{equation}\label{LCLTeq2}
\pr{S(2n)=2k} \leq \frac{K}{\sqrt{n}}\exp\left(-\frac{k^2}{n}\right).
\end{equation}
\end{thm}

\begin{proof}
$\pr{S(2n)=\pm 2n} = \left(\frac{1}{2}\right)^{2n}$ and for $|k|<n$, Stirling's formula yields
\begin{eqnarray}\label{LCLTprob}
\pr{S(2n)=2k} & = & {2n \choose
  n+k}\left(\frac{1}{2}\right)^{2n}\nonumber \\
& = & \frac{\sqrt{4\pi
    n}}{\sqrt{2(n+k)\pi}\sqrt{2(n-k)\pi}}\frac{(2n)^{2n}}{(n+k)^{n+k}(n-k)^{n-k}}\left(\frac{1}{2}\right)^{2n}\left(1+\bigo{\frac{1}{n}}\right)\nonumber \\
& = & \frac{\sqrt{\pi
    n}}{\sqrt{(n+k)\pi}\sqrt{(n-k)\pi}} \phi(n,k)\left(1+\bigo{\frac{1}{n}}\right),
\end{eqnarray}
where $\phi(n,k) = \frac{n^{2n}}{(n+k)^{n+k}(n-k)^{n-k}}$. The Taylor
expansion of $\ln(1+x)$ gives for $|k| < n$

\begin{eqnarray*}
\log \phi(n,k) & = & 2n\log n - (n+k)\log(n+k) - (n-k)\log(n-k)\\
& = & -n(\log(1+\frac{k}{n}) + \log(1-\frac{k}{n})) - k(\log(1+\frac{k}{n})
- \log(1-\frac{k}{n}))\\
& = & -n\left(\S{l=1}{\infty}{\frac{(-1)^{l+1}}{l}\left(\frac{k}{n}\right)^l} -
\S{l=1}{\infty}{\frac{1}{l}\left(\frac{k}{n}\right)^l}\right)\\
& & \hspace{4pc}-k\left(\S{l=1}{\infty}{\frac{(-1)^{l+1}}{l}\left(\frac{k}{n}\right)^l}
+ \S{l=1}{\infty}{\frac{1}{l}\left(\frac{k}{n}\right)^l}\right)\\
& = & 2\left(n\S{l=1}{\infty}{\frac{1}{2l}\left(\frac{k}{n}\right)^{2l}} -
k\S{l=1}{\infty}{\frac{1}{2l-1}\left(\frac{k}{n}\right)^{2l-1}}\right) = -\S{l=1}{\infty}{\frac{1}{l(2l-1)}\frac{k^{2l}}{n^{2l-1}}},
\end{eqnarray*}
which implies that 
\begin{equation}\label{phi(n,k)}
\phi(n,k) =
\exp\left\{-\S{l=1}{\infty}{\frac{1}{l(2l-1)}\frac{k^{2l}}{n^{2l-1}}}\right\}.
\end{equation}
In particular, there exists a $K>0$ such that for all $n\in\N, N\geq 2$,
and $|k|\leq n^{(2N-1)/2N}$, 
\begin{equation}\label{part1}
\phi(n,k) =
\exp\left\{-\S{l=1}{N-1}{\frac{1}{l(2l-1)}\frac{k^{2l}}{n^{2l-1}}}\right\}
(1+\eta_N(n,k)),
\end{equation}
where $\eta_N(n,k) \leq K\frac{k^{2N}}{n^{2N-1}}$.
% = o(1)$.
Also, 
\begin{eqnarray*}
\frac{\sqrt{\pi n}}{\sqrt{(n+k)\pi}\sqrt{(n-k)\pi}} & = &
 \sqrt{\frac{1}{\pi}}\sqrt{\frac{1}{n(1-k^2/n^2)}}\\
& = & \sqrt{\frac{1}{\pi
n}}\left(1+\S{l=1}{\infty}{\left(\frac{k}{n}\right)^{2l}\frac{(2l-1)(2l-3)\cdots
    3\cdot 1}{2^l l!}}\right),
\end{eqnarray*}
which converges for $|k|<n$, so for any $b<1$, there is a constant $C_b$
such that for all $|k|\leq bn$, 
$$\left|\frac{\sqrt{\pi n}}{\sqrt{(n+k)\pi}\sqrt{(n-k)\pi}} - \sqrt{\frac{1}{\pi
n}}\right| \leq C_b \frac{k^2}{n^{5/2}}.$$
In particular, there exist a constant $K>0$ and a function $\eps(x)$ with
$\eps(x) \leq Kx$ such that for every $a <
1$, every $n\geq 1$, and every $|k|\leq n^a$, 
\begin{equation}\label{part2}
\frac{\sqrt{\pi n}}{\sqrt{(n+k)\pi}\sqrt{(n-k)\pi}} = \sqrt{\frac{1}{\pi
n}}\left(1+\eps\left(\frac{k^2}{n^2}\right)\right) = \sqrt{\frac{1}{\pi
n}}(1+o(1)).
\end{equation}
Moreover, it is easy to see that for all $|k|<n$.
\begin{equation}\label{trivialbound}
\frac{\sqrt{\pi n}}{\sqrt{(n+k)\pi}\sqrt{(n-k)\pi}} \leq 1.
\end{equation}

\eqref{LCLTprob}, \eqref{part1}, and \eqref{part2} yield \eqref{LCLTeq1}. To
show \eqref{LCLTeq2}, it suffices to consider separately two cases: If
$k\leq n^{5/6}$, \eqref{part1} and \eqref{part2} directly yield the
result. If $n^{5/6}\leq k \leq n$, \eqref{phi(n,k)} and
\eqref{trivialbound} show that $$\pr{S(2n)=2k} \leq
K\exp(-\frac{k^2}{n})\exp(-\frac{k^4}{6n^3}) \leq \frac{K}{\sqrt{n}}\exp(-\frac{k^2}{n}).$$

\end{proof}

\begin{remark} Replacing $k$ by $n$ in \eqref{phi(n,k)} yields, as one
  would expect, $\phi(n,n) = \left(\frac{1}{2}\right)^{2n}$ 
\end{remark}

\begin{cor}[Large deviations for one-dimensional simple random walk]\label{largedevsrw}
If $S$ is one-dimensional simple random walk, there exists a $K>0$
such that for every $n, k\in\N$, every $r\geq 1$, 
$$\pr{|S(2n)|\geq r\sqrt{n}}\leq K\exp\left\{-\frac{r^2}{4}\right\}.$$
\end{cor}

\begin{proof}
Given Theorem \ref{LCLTthm}, all that is left to do is integrate: 
\begin{eqnarray*}
\pr{|S(2n)|\geq r\sqrt{n}} & \leq &
2\S{k=[r\sqrt{n}/2]}{n}{\pr{S(2n)=2k}}\\
& \leq &
\frac{K}{\sqrt{n}}\S{k=[r\sqrt{n}/2]}{n}{\exp\left\{-\frac{k^2}{n}\right\}} \leq
K\int_{r/\sqrt{2}}^{2\sqrt{2n}}e^{-x^2/2} \,dx \leq Ke^{-r^2/4}.
\end{eqnarray*}
\end{proof}

\begin{remark} The corollary is equivalent to the statement that there exists a $K>0$
such that for every $n, k\in\N$, every $r\geq 1$, 
$$\pr{|S(n)|\geq r\sqrt{n}}\leq K\exp\left\{-\frac{r^2}{2}\right\},$$
which is in agreement with the bound for one-dimensional Brownian motion.
\end{remark}

We now give the analogue of Lemma \ref{bmtoofar} for random walk.
\begin{lem}\label{rwtoofar} There exists a positive constant $K$ such
  that if $S$ is planar simple random walk, then 
\begin{eqnarray*}
& (a) & \pr{|S(2n)|\geq r\sqrt{n}} \leq Ke^{-r^2/4}.\\
& (b) & \pr{\sup_{0\leq t \leq n}|S(2t)| \geq r\sqrt{n}} \leq 2Ke^{-r^2/4}.
\end{eqnarray*}
\end{lem}

\begin{proof}
We know from Section \ref{1-2} that we can write $S(n) =
(S^1(n),S^2(n))$, where $S^1(n)$ and $S^2(n)$  are
independent random walks on
$\frac{1}{\sqrt{2}}e^{i\pi/4}\cdot\Z$ and
$\frac{1}{\sqrt{2}}e^{i3\pi/4}\cdot\Z$, respectively. These are random
walks on shrunken lattices and we get from Corollary \ref{largedevsrw}
that for $k=1, 2$, 
$$\pr{|S^k(2n)|\geq r\sqrt{\frac{n}{2}}}\leq
K\exp\left\{-\frac{r^2}{4}\right\},$$ 
which, together with the obvious inequality $\pr{|S(2n)|\geq r\sqrt{n}}
\leq 2\pr{|S^1(2n)| \geq r\sqrt{n/2}}$, gives part $(a)$. Part $(b)$
follows from Lemma \ref{reflection}.

\end{proof}

\section{Beurling estimates}\label{beurling}

It is often useful to know how likely it is for Brownian motion to get to distance $R$ without hitting a set $A$ with $d(B(0),A)=r\leq R$ and rad$(A)\geq 2R$. The probability of this event can be bounded above by a power function of the ratio $\frac{r}{R}$, uniformly for all sets $A$. Given the Beurling Projection Theorem which we state below, it is easy to find the best possible exponent of this power function and we do it in this section. The discrete case is more difficult and we just refer the reader to \cite{kesten}, where the proof is given.

The first result of this section is the Beurling Projection Theorem (for a proof, see \cite{bass} or \cite{ahlfors2}). It says that among all connected sets of a given radius, that which Brownian motion will most likely avoid is a straight line. Let $\D$ be the closed unit disk centered at the origin and consider a set $E \subset R\D = \{z\in\C:|z|\leq R\}$. The circular projection of $E$ is defined to be $\gam(E)=\{|z|:z\in E\}$. For a set $A\subset \C$, let $T_A = \inf\{t\geq 0:B(t)\in A\}$ and $\Xi_R = T_{\bd R\D}$.
\begin{thm}[Beurling Projection Theorem]\label{projection}
For all $R\geq 1$ and $E\subset \C$,
$$\pro{-1}{\Xi_R<T_E} \leq \pro{-1}{\Xi_R<T_{\gam(E)}}.$$
\end{thm}

We will be interested in the case where $E$ satisfies $\gam(E)=[0,R]$. Now that we have Theorem \ref{projection}, we know that finding an upper bound for $\pro{-1}{\Xi_R<T_{[0,R]}}$ also provides an upper bound for $\pro{-1}{\Xi_R<T_E}$ for all sets $E \subset R\D$ with $\gam(E)=[0,R]$. We can compute such a bound via a sequence of conformal maps, using the fact that the exit distribution of the upper half-plane is a Cauchy distribution, and the fact that harmonic measure is conformally invariant (see \cite{bass} for a proof of this). 
%It turns out that the bound we find is optimal up to a multiplicative constant.

\begin{figure}
 \psfragscanon
         \centering
         \includegraphics[height=2in,width=2.6in]{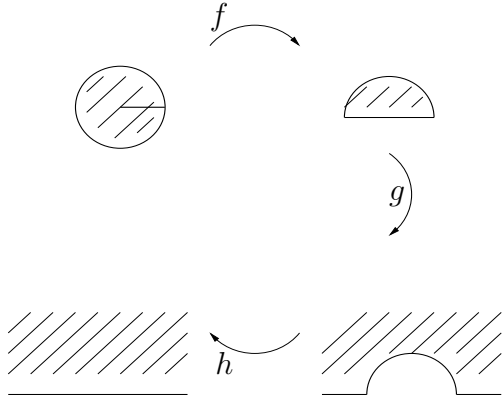}
         \caption{The sequence of conformal transformations leading to the Beurling estimate.}
         \label{Fig4}
\end{figure}

Consider the following domains, where $\opdisk = \{z\in\C:|z|<1\}$ is the open unit disk and recall that $\D=\bd\opdisk$:
the upper half-plane $\h = \{z\in\C:\im{z}>0\}$, the slit unit disk $\opdisk_s = \opdisk \setminus \{z\in\C:0\leq \re{z} < 1 ; \im{z}=0\}$, the upper half-disk $\opdisk_u = \opdisk \cap \h$, and the complement in $\h$ of the closed upper half-disk $\h_\opdisk = \h \cap \{z\in\C:|z|>1\}.$ These domains are linked by the following conformal transformations (surjective conformal maps):
$$\opdisk_s \stackrel{f}{\rightarrow} \opdisk_u \stackrel{g}{\rightarrow} \h_\opdisk \stackrel{h}{\rightarrow} \h,$$
where $f(z) = \sqrt{z}, g(z) = -\frac{1}{z}, h(z)=z+\frac{1}{z}$. Then $h \circ g \circ f (-\eps) = (\frac{1}{\sqrt{\eps}}-\sqrt{\eps})i$ and  $h \circ g \circ f (\opdisk) = [-2,2].$ See Figure \ref{Fig4}.

Conformal invariance of Brownian motion implies that 
$$\pro{\eps}{B(T_{\partial \opdisk_s})\in \D} = \pro{(\frac{1}{\sqrt{\eps}}-\sqrt{\eps})i}{B(T_{\R}) \in [-2,2]}.$$

Using the fact that the exit distribution of the upper half-plane is a Cauchy distribution and the Beurling Projection Theorem, Brownian scaling gives the following:
\begin{thm}[Continuous Beurling Estimate] \label{beurlingcontinuous} There exists a constant $K >0$ such that for any $R \geq 1$, any $x$ with $|x|\leq R$, any set $A$ with $[0,R] \subset \gam(A)$,
$$\pro{x}{\Xi_R \leq T_A}\leq K\left(\frac{|x|}{R}\right)^{1/2}.$$
\end{thm}

As we pointed out earlier, extending this result to the random walk case is not as easy, mainly because none of the conformal invariance techniques are available. In \cite{kesten}, Kesten first showed that the Beurling estimate holds in the discrete case as well. We state the theorem here without a proof.
For $R\in \R_+$, we define ${\mathcal A}_R$ to be the set of subsets $A$ of $\Z^2$ containing 0 and for which $\sup\{|x|:x\in A\}\geq R$. We also let for $A\subset\Z^2, \tau_A = \inf\{k\geq 0:S(k)\in A\}$ and $\xi_R = \inf\{k\geq 0:|S(k)|>R\} $.

\begin{thm}[Discrete Beurling Estimate]\label{beurlingdiscrete}
 Let $\tau_A = \inf\{n\geq 1:S(n) \in A\}$ and $\xi_R = \inf\{k\geq 0:|S(k)|\geq R\}$. Then there exists a constant $K > 0$ such that if $A \in {\mathcal A}_R, |x|<R,$
$$\pro{x}{\xi_R\leq\tau_A}\leq K\left(\frac{|x|}{R}\right)^{1/2}.$$
\end{thm}

Note that although the exponents are the same in the continuous and discrete case, it is not clear that in a given discrete disk, a straight line is the easiest set to avoid for random walk.

%\begin{thm}[Discrete Beurling Estimate]\label{beurlingdiscrete}
% Let $R\in \R$ be positive and define ${\mathcal A}_R$ to be the set of subsets of $\Z^2$ for which $\sup\{|x|:x\in A\}=R$. Let $\tau_A = \inf\{n\geq 1:S(n) \in A\}$ and $\xi_R = \inf\{k\geq 0:|S_k|\geq R\}$. Then there exists a constant $K > 0$ such that if $A \in {\mathcal A}_R, |x|<R,$
%$$\pro{x}{\tau_A>\xi_R}\leq K\left(\frac{|x|}{R}\right)^{1/2}.$$
%\end{thm}

%Note that although the exponents are the same in the continuous and discrete case, it is not clear that in a given discrete disk, a straight line is the easiest set to avoid for random walk.

\section{Skorokhod Embedding}\label{Skorokhod}  % use *-form to suppress numbering

%Since the 1960's, much has been written about ways to couple random walk with Brownian motion and the subject is now quite well understood. Usually the goal is to put the two processes on a same probability space in such a way that with large probability they remain close to each other at all times of a given time interval. For an extensive discussion of this problem, see \cite{approx}. Two such couplings of somewhat different nature will be needed in this thesis and we discuss them in this chapter. The first is based on the Skorokhod embedding scheme. The second is much sharper and is the so-called KMT approximation.

Knowledge about the behavior of Brownian motion can often be used to
understand the behavior of random walk, and vice-versa. Coupling
arguments turn out to be particularly useful in many cases. We briefly
discuss one of them here, namely Skorokhod embedding.

\subsection{The one-dimensional case}

%We state here Skorokhod's theorem on how to embed general random walk in Brownian motion but will not need it directly since we will make an explicit construction of the coupling. Note that in this section, random walk and Brownian motion are one-dimensional.

%\begin{thm}[Skorokhod Embedding]\label{sko} Suppose that $(X_i)_{i\geq 1}$ are independent, identically distributed real-valued random variable with mean 0 and variance 1. Then there exist a probability space containing a Brownian motion $\{B(t):t\geq 0\}$, the random variables $(X_i)_{i\geq 1}$, a sequence of stopping times $0=T_0 \leq T_1 \leq T_2 \leq ...$, such that the increments $T_n-T_{n-1}$ are independent, identically distributed, $\ex{T_n} = n$, and the sequence $\{B(T_n)\}_{n\geq 1}$ has the same distribution as the random walk $\{S(n)\}_{n\geq 1}$ associated with $\{X_i\}_{i\geq 1}$.
%\end{thm}

%For more details, see for instance \cite{durrett}. Note that although this theorem is a beautiful result, it does not say anything about how close the random walk and the Brownian motion actually are if we look at them path by path. For theorems addressing this question, we again refer the reader to \cite{approx}.

\begin{thm}\label{skocoupling} There exists a probability space containing a linear standard
  Brownian motion $B$ and a one-dimensional simple random walk $S$ such that for every $g(n)\geq 1$,
  %satisfying $g(n) = \bigo{n^{1/4}}$
  there exist constants $b,K>0$ such that 
$$\pr{\sup_{0\leq t \leq n}|B(t) - S(t)|\geq n^{1/4}g(n)} \leq Kne^{-bg(n)}.$$
\end{thm}

\begin{proof} 
We define $T_0 = 0$ and for $i\geq 1, T_i = \inf\{t\geq T_{i-1}:|B(t) -
B(T_{i-1})| = 1\},$ and define the simple random walk $S(n) :=
B(T_n)$. For notational purposes, we let $h(n) = g(n)\sqrt{n}$ and for $k\geq 1$, 
$$I_k=\left[(k-1)[h(n)],k[h(n)]\right].$$ For $n\geq 4$, we have the covering $[0,n] \subset \displaystyle{\bigcup_{k=1}^{[\sqrt{n}/g(n)]+3}I_k}$. Also, if $0\leq t-s \leq h(n)$, the interval $[s,t]$ intersects at most 3 of the $I_k$'s. We use this and Corollary \ref{bmstoppingtimes} to see that there exist constants $K, a >0$ such that
\begin{align*}
\p\{\sup_{0\leq k \leq n} & |B(k) - S(k)|\geq \frac{1}{2}n^{1/4}g(n)\}\\
& \leq  \pr{\max_{1\leq k\leq n}|T_k-k|\geq h(n)} +
\pr{\mathop{\sup_{0\leq s\leq n}}_{|t-s|\leq h(n)}|B(t)-B(s)|\geq \frac{1}{2}n^{1/4}g(n)}\\
& \leq  Kne^{-ag(n)} + \pr{\mathop{\sup_{1\leq k\leq [\sqrt{n}/g(n)]+3}}_{t\in I_k}|B(t)-B((k-1)[h(n)])|\geq \frac{1}{6}n^{1/4}g(n)}\\
& \leq  K\left(ne^{-ag(n)} + \frac{\sqrt{n}}{g(n)}\pr{\sup_{0\leq t\leq [h(n)]}|B(t)|\geq \frac{1}{6}n^{1/4}g(n)}\right)\\
& \leq  K\left(ne^{-ag(n)} + \frac{\sqrt{n}}{g(n)^{3/2}}e^{g(n)/72}\right) \leq Kn e^{-bg(n)},
\end{align*}
where the two last inequalities follow from Brownian scaling, Lemma \ref{bmtoofar}, and the choice $b=\min\{a,1/72\}$.
To get the result for all real times $t\in [0,n]$, it suffices to observe that
\begin{align*}
\mathbb{P}\{&\sup_{0\leq t \leq n}|B(t) - S(t)|\geq n^{1/4}g(n)\}\\
& \leq \pr{\displaystyle{\sup_{0\leq k \leq n}|B(k) - S(k)|\geq \frac{1}{2}n^{1/4}g(n)}}\\
& \hspace{5pc}+ n\pr{\sup_{0\leq t\leq 1}|B(t)-S(t)|\geq \frac{1}{2}n^{1/4}g(n)}\\
& \leq K\left(ne^{-bg(n)} + ne^{-b'n^{1/2}}\right) \leq Kne^{-bg(n)}.
\end{align*}
\end{proof}

%\begin{remark} The theorem implies that for any $\phi$ such that $\phi(n) \to \infty$ as $n\to\infty$, 
%$$\sup_{0\leq t \leq n}\frac{|B(t)-S(t)|}{n^{1/4}\log n\phi(n)}\stackrel{P}{\to}0.$$
%This is close to being the best we can hope for when using Skorokhod's method. In \cite{kiefer}, Kiefer shows that the best possible result for simple random walk with Skorokhod's method is
%$$\frac{|S(n)-B_n|}{(n\log\log n)^{1/4}(\log n)^{1/2}}\stackrel{\text{a.s.}}{=}\bigo{1}.$$
%\end{remark}

\subsection{Extending the result to the plane}

The construction we made in the previous section does not work in dimensions other than 1. However, there is a simple way of getting around this 
problem, which is as follows:

Let $B^1(t), B^2(t)$ be independent one-dimensional Brownian motions. Then $$B(t)=(B^1(t),B^2(t))$$ is a planar Brownian motion. For $i\geq 0, j=1,2$, let $T_i^j$ be the stopping times for $B^j(t)$ as defined in the previous subsection and for $i\geq 0$, define $S_n^j = B^j(T_n^j)$. Then $S^1, S^2$ are independent one-dimensional random walks. We let $L_i=(L_i^1,L_i^2)$ be independent random vectors, independent of $B^1,B^2$, with distribution
$$\pr{L_i=(1,0)}=\pr{L_i=(0,1)}=\frac{1}{2}.$$
If we define $U_n^j = \S{i=1}{n}{L_i^j}$, it is easy to check that $S(n):=(S^1(U_n^1),S^2(U_n^2))$ is a planar simple random walk. The statement and the proof of the main result are essentially the same as in one dimension. The only difference is that now, the time-parameter is different for the two processes. A heuristic reason for the different time scales is that Brownian motion moves a little bit faster since it is not restricted to moving along the lines of the lattice, but can take ``diagonal shortcuts''. See also Lemma \ref{meansquaredistance}.

\begin{thm}
 There exists a coupling of standard Brownian motion $B$ and simple
 random walk $S$ in the plane such that for all $g(n)\geq 1$ satisfying $g(n) = \bigo{n^{1/4}}$, there exist constants $b,K>0$ such that 
$$\pr{\sup_{0\leq t \leq n}|B(t) - S(2t)|\geq n^{1/4}g(n)} \leq Kne^{-bg(n)}.$$
\end{thm}

\begin{proof}
If $S$ is defined as above, we have
$$\pr{\max_{1\leq k \leq n}|B(k) - S(2k)|\geq n^{1/4}g(n)} \leq 2\pr{\max_{1\leq k \leq n}|B^1(k) - S^1(T_{2k}^1)|\geq \frac{1}{\sqrt{2}}n^{1/4}g(n)}.$$

Since $T_{2k}$ is a sum of $2k$ random variables of mean $1/2$ and finite variance, the exact same argument as in Theorem \ref{skocoupling} can be used to conclude the proof.
\end{proof}

\section{Exiting distributions for simple random walk}

The Laplacian $\Delta$ of a $C^2$ function $f:\R^2\to\R$ is defined by
$$\Delta f(x,y) = \frac{\bd{}^2f}{\bd{x}^2} + \frac{\bd{}^2f}{\bd{y}^2}(x,y).$$
If in a domain $D, \Delta f \equiv 0$, we say that $f$ is harmonic in $D$. Similarly, for any function $f:\Z^2 \to \R$ the discrete Laplacian is
$$\Delta f(x,y) = \frac{1}{4}\sum (f(x',y')-f(x,y)),$$
where the sum is over $\{(x',y'):|(x',y')-(x,y)|=1\}$. By analogy with the continuous case, we say that $f:\bar{D}\to\R$ is discrete harmonic in a set $D$ if $\Delta f \equiv 0$ in $D$.

To solve the Dirichlet problem in a domain $D$ with boundary condition $\phi$, where $\phi:\bd{D}\to\R$ is to find a function $u:\bar{D}\to \R$ such that $u$ is harmonic in $D$ and $u\equiv \phi$ on $\bd{D}$. The discrete Dirichlet problem is defined in the natural analogous way.

The Dirichlet problem is intimately related to Brownian motion and so is the discrete Dirichlet problem to random walk. Solving the Dirichlet problem with appropriate boundary conditions is equivalent to computing the probability that Brownian motion (or random walk in the discrete case) leaves a domain at a given subset of the boundary. More precisely, if $A$ and $B$ are disjoint subsets of the boundary of $D$ and $A\cup B = \bd{D}$, solving the Dirichlet problem with boundary value 1 on $A$ and 0 on $B$ is equivalent, under some mild assumptions, to finding the probability that Brownian motion (in the continuous case) or random walk (in the discrete case) leaves $D$ at $A$. (See \cite{durrett2} for a discussion of this in the continuous case and \cite{greenbook} for the discrete case.) 
%We will need upper bounds for various such exiting probabilities in this thesis and solve the corresponding Dirichlet problems here. We also prove a ``difference estimate'' which gives a bound on the difference between the exiting distribution of a domain when starting at neighboring points. While working on the problems of this thesis, we initially believed that we would need this estimate, which eventually was replaced by another. As the result is of interest by itself and is not, to our knowledge, anywhere in the literature, we present it here as well.

%\section{Dirichlet Problem in the (Infinite) Rectangle}

%Results from \cite{sle}.

\subsection{Discrete Dirichlet problem in the finite and infinite rectangles}

We will first solve the discrete Dirichlet problem on a discrete rectangle with boundary conditions a general function $\phi$ on one side and 0 on the others. We will then use this to find bounds for the problem with specific $\phi$ and points inside the rectangle.

We start with a trivial lemma which will be needed in our study of the discrete Dirichlet problem in some specific domains.

\begin{lem}\label{a_jlemma} If for $1 \leq j \leq n-1$, $a_j = a_j(n)$ is defined to be the positive solution of the equation
\begin{equation}\label{110306}
\cosh(a_j)=2-\cos(\frac{\pi j}{n}),
\end{equation}
then 
%$$a_j = \frac{\pi j}{n}(1 + \bigo{\frac{j}{n}}).$$
%Moreover, 
for any $1 \leq j \leq n-1$,
$$\frac{j}{2n} \leq a_j \leq \frac{\pi j}{n}.$$
\end{lem}

\begin{proof}
The equality
$$\displaystyle{\sum_{k\geq 1}\frac{a_j^{2k}}{(2k)!} = \sum_{k\geq 1}(-1)^{k+1}\frac{(\pi j/n)^{2k}}{(2k)!}},$$
obtained by Taylor-expanding \eqref{110306}, allows us to see that $\forall \; n>0, \forall \; j \leq n/\pi,$ 
$$\frac{\pi j}{2n} \leq a_j \leq \frac{\pi j}{n}.$$
Indeed, suppose that for some $n>0$ and $j\leq n\pi, a_j > \frac{\pi j}{n}$. Then 
$$\displaystyle{\sum_{k\geq 1}\frac{(\pi j/n)^{2k}}{(2k)!} < \sum_{k\geq 1}(-1)^{k+1}\frac{(\pi j/n)^{2k}}{(2k)!}}.$$
This is clearly impossible. Also, if we suppose that for some $0\leq a_j \leq \frac{\pi j}{2n}$, we get the inequality
$$0 \leq \S{k\geq 1}{}{\frac{(\pi j/n)^{2k}\left((1/2)^{2k}-(-1)^{k+1}\right)}{(2k)!}}.$$
The first term in this sum is $\frac{-3}{8}(\frac{\pi j}{n})^2$. The sum of the positive terms is 
$$\leq \frac{17}{16}(\frac{\pi j}{n})^4\S{k\geq 1}{}{\frac{1}{(4k)!}} \leq \frac{1}{10}(\frac{\pi j}{n})^4,$$ 
and we get a contradiction. It is also easy to see directly from the $\cos$ and $\cosh$ functions that for $n/\pi \leq j \leq n-1, \, 1/2 \leq a_j \leq \pi$. This is not optimal but sufficient for our needs. The lemma now follows easily.

%gives the following lower bound for all $j\leq n-1$:

%\begin{equation}\label{a_jlowerbound}
%a_j \geq \frac{j}{2n}.
%\end{equation}
\end{proof}

We let 
$$R(l,n)=\{(x,y) \in \Z^2:1\leq x\leq l-1,1\leq y \leq n-1\},$$ 
be the discrete rectangle of ``side lengths'' $l-1$ and $n-1$, with boundary $\partial R(l,n)$ and closure $\bar{R}(l,n)=R(l,n) \cup \partial R(l,n)$.
\begin{lem}\label{Rectangle} Let $\phi:\{1,..,n-1\} \to \R$ be a given function. Then the unique function $f:\bar{R}(l,n) \to \R$ satisfying
$$\Delta f(x,y) = 0 \text{ in } R(l,n),$$
$$f(x,y) = \left\{   
\begin{array}{ll}
                \phi(y) & \text{ on } \{(l,y):1\leq y \leq n-1\}\\
                0 & \text{ on } \partial R(l,n) \setminus \{(l,y):1\leq y \leq n-1\}
\end{array} \right.,$$
is given by 
\begin{equation}\label{rectangle}
f(x,y) = \displaystyle{\sum_{j=1}^{n-1}b_j(\phi)\frac{\sinh(a_jx)}{\sinh(a_jl)}\sin(\frac{\pi yj}{n})},
\end{equation}
where $a_j$ is the positive solution of 
\begin{equation}\label{a_j}
\cosh(a_j) = 2 - \cos(\frac{\pi j}{n}),
\end{equation}
and 
\begin{equation}\label{b_j}
b_j(\phi)=\frac{2}{n-1}\displaystyle{\sum_{y=1}^{n-1}\phi(y)\sin(\frac{\pi yj}{n})}.
\end{equation}
\end{lem}

\begin{proof}
It suffices to check that the given function is harmonic and that it has the right values on $\partial R(l,n)$. Uniqueness follows from \cite[Theorem 1.4.5]{greenbook}.

We first check the boundary conditions. It is clear that 
$$f(0,y) = f(x,0) = f(x,n) = 0 \;\; \forall \; x,y \in \Z.$$
Also,
\begin{eqnarray*}
f(l,y) & = & \displaystyle{\sum_{j=1}^{n-1}b_j(\phi)\sin(\frac{\pi yj}{n})}\\
& = & \frac{2}{n-1}\displaystyle{\sum_{j=1}^{n-1}\sum_{y=1}^{n-1}\phi(y)\sin^2(\frac{\pi yj}{n})}.
\end{eqnarray*}
It is easy to see that if $1 \leq y \leq n-1$, then
$$\displaystyle{\sum_{j=1}^{n-1}\sin^2(\frac{\pi yj}{n})} = \frac{n-1}{2},$$
from which it follows that 
$$f(l,y) = \phi(y) \;\; \forall \; y \in \{1,..,n-1\}.$$
To see that $f$ is harmonic in $R(l,n)$, we do a straightforward computation: Fix $(x,y) \in R(l,n)$. Then, using the fact that $\sin(x+y)=\sin(x)\cos(y)+\sin(y)\cos(x)$,
\begin{align*}
\Delta &f(x,y)\\
& = f(x+1,y) + f(x-1,y) + f(x,y+1) + f(x,y-1) - 4f(x,y)\\
& = \displaystyle{\sum_{j=1}^{n-1}b_j\left[2\sin(\frac{\pi yj}{n})(\cos(\frac{\pi j}{n})-1)\frac{\sinh(a_jx)}{\sinh(a_jl)} + 2\sin(\frac{\pi yj}{n})\frac{\sinh(a_jx)}{\sinh(a_jl)}(\cosh(a_j)-1)\right]}\\
& = 2\displaystyle{\sum_{j=1}^{n-1}b_j\left[\sin(\frac{\pi yj}{n})\frac{\sinh(a_jx)}{\sinh(a_jl)}\left((\cos(\frac{\pi j}{n}) - 1)+(\cosh(a_j)-1)\right)\right]}.
\end{align*}
We see that this is 0 if $\cos(\frac{\pi j}{n})+\cosh(a_j)=2$.

\end{proof}

We now turn to the infinite rectangle
$$\mathcal{R}(n)= \{(x,y) \in \Z^2:x \geq 1, 1 \leq y \leq n-1\}$$
with boundary $\partial\mathcal{R}(n)$ and closure $\bar{\mathcal{R}}(n) = \mathcal{R}(n) \cup \partial\mathcal{R}(n)$.

\begin{lem}\label{Infrectangle}Let $\phi:\{1,..,n-1\} \to \R$ be a given function. Then the unique bounded function $f(x,y):\bar{\mathcal{R}}(n) \to \R$ satisfying
$$\Delta f(x,y) = 0 \text{ in }\mathcal{R}(n),$$
$$f(x,y) = \left\{   
\begin{array}{ll}
                \phi(y) & \text{ on } \{(0,y):1\leq y \leq n-1\}\\
                0 & \text{ on } \partial \mathcal{R}(n) \setminus \{(0,y):1\leq y \leq n-1\}
\end{array} \right.$$
is given by 
\begin{equation}\label{infrectangle}
f(x,y) = \displaystyle{\sum_{j=1}^{n-1}b_j(\phi)\exp(-a_jx)\sin(\frac{\pi yj}{n})},
\end{equation}
where $a_j$ is the positive solution of 
\begin{equation*}
\cosh(a_j) = 2 - \cos(\frac{\pi j}{n}),
\end{equation*}
and $$b_j(\phi)=\frac{2}{n-1}\displaystyle{\sum_{y=1}^{n-1}\phi(y)\sin(\frac{\pi yj}{n})}.$$
\end{lem}

\begin{proof}

We invoke \cite[Theorem 1.4.8]{greenbook} to show uniqueness. We can check that $f$ has the right boundary conditions exactly as in Lemma \ref{Rectangle}. Harmonicity follows from 
$$\Delta f(x,y) = \displaystyle{\sum_{j=1}^{n-1}b_j(\phi)\exp(-a_jx)\sin(\frac{\pi yj}{n})\left(2\cos(\frac{\pi j}{n})+2\cosh(a_j)-4\right)}.$$
\end{proof}

We now find upper bounds for solutions of the Dirichlet problem in a finite and infinite rectangle at particular points in the case where $\phi \equiv 1$.

\begin{lem}\label{f(1,y)}If $f(x,y)$ is the solution of the Dirichlet problem in $R([an],n)$ with $\phi(y) \equiv 1$, then there exists a positive constant $K$, depending on $a$, such that for all $y$ and all $n,$
$$f(1,y) \leq K\frac{y}{n^2}.$$
\end{lem}

\begin{proof}
This is a particular case of Lemma \ref{Rectangle}. First note that $a_j$ and $b_j$ depend on $n$.
$$b_j:=b_j(\phi) = \frac{2}{n-1}\displaystyle{\sum_{k=1}^{n-1}\sin \left(\frac{\pi kj}{n}\right) \stackrel{n \to \infty}{\longrightarrow} 2\int_0^1 \sin(\pi jx)\,dx} = 
\left\{ 
\begin{array}{ll} 
\frac{4}{\pi j}, & j \text{ odd }\\
0, & j \text{ even, }
\end{array}, \right.$$
so $\exists \, C>0, \text{ s.t. }\forall \; n>0, \; \forall \; j \in \{1,..,n-1\}$, 
\begin{equation}\label{b_jbound}
b_j \leq \frac{C}{j}.
\end{equation}
We also note that 
\begin{equation}\label{sin}
\sin(x) \leq x, \; \forall x \geq 0.
\end{equation}
Finding a bound for the term $\frac{\sinh(a_jx)}{\sinh(a_j[an])}$ is more delicate. We first note that for all $x \geq 0, \, \sinh(x) \geq x$, and for $x\leq 1,\; \sinh(x)\leq 2x.$ We also recall from Lemma \ref{a_jlemma} that $\frac{j}{2n} \leq a_j \leq \frac{\pi j}{2n}.$

%The equality

%$$\displaystyle{\sum_{k\geq 1}\frac{a_j^{2k}}{(2k)!} = \sum_{k\geq 1}(-1)^{k+1}\frac{(\pi j/n)^{2k}}{(2k)!}},$$

%By supposing that for some $j$ and $n, a_j \geq \frac{\pi j}{n}$ and arriving to a contradiction, 

%obtained from (\ref{a_j}), allows us to see (e.g. by contradiction) that $\forall \; n>0, \forall \; j \leq n/\pi,$ 

%$$\frac{\pi j}{2n} \leq a_j \leq \frac{\pi j}{n}.$$

%It is also easy to see that for $n/\pi \leq j \leq n-1, a_j \geq 1/2$.

%This gives the following lower bound for all $j\leq n-1$:

%where the lower bound is found in the same way as the upper bound.
\smallskip

\begin{itemize}
\item If $1 \leq j < \frac{8}{a\pi}$, then $a_j[an] \geq C$ for some $C>0$ depending on $a$, and $a_j$ is small, so that
$$\frac{\sinh(a_j)}{\sinh(a_j[an])} \leq \frac{C_1a_j}{C_2} \leq C\frac{j}{n}.$$

\item If $\frac{8}{a\pi} \leq j < \frac{n}{\pi}$ and $n$ is large enough,
$$a_j[an] \geq [an]\frac{j}{2n} \geq cj \geq 1,$$
for some $c>0$ independent of $j$, so that, since $a_j < 1,$
$$\frac{\sinh(a_j)}{\sinh(a_j[an])} \leq C\frac{j/n}{\exp(cj)}.$$

\item Finally, if $n/\pi \leq j \leq n-1$, then
$$\frac{\sinh(a_j)}{\sinh(a_j[an])} \leq C_1\exp(-C_2n).$$
\end{itemize}

%\medskip

%For $\frac{n}{\pi[an]} < j < \frac{n}{\pi}, a_j < 1$ and $a_j [an]>2$, so that 

%$$\frac{\sinh(a_j)}{\sinh(a_jL)} \leq \frac{a_j + \bigo{a_j^3}}{a_jL} = \frac{\sqrt{2}\pi j/n+\bigo{(j/n)^3}}{\sqrt{2}\pi jL/n+\bigo{(jL/n)^3}},$$

%and letting $L = [an]$, we find that 

%$$\frac{\sinh(a_j)}{\sinh(a_jL)} \leq C\frac{j/n + \bigo{(j/n)^3}}{j(1+\bigo{j^2})} \leq C\left(\frac{1}{nj^2}+\frac{1}{n^3}\right).$$
Plugging these bounds, as well as (\ref{b_jbound}) and (\ref{sin}) into (\ref{rectangle}), we get
%$$f(1,y) \leq \displaystyle{\sum_{j=1}^{n-1}\frac{1}{j}\frac{1}{nj^2}\left(\frac{yj}{n}+\bigo{(\frac{j}{n})^2}\right)} \leq \frac{Cy}{n^2}$$
\begin{eqnarray*}
f(x,y) & = & \displaystyle{\sum_{j=1}^{n-1}b_j(\phi)\frac{\sinh(a_jx)}{\sinh(a_j[an])}\sin(\frac{\pi yj}{n})}\\
& \leq & C\left(\displaystyle{\sum_{j=1}^{[4/a\pi]}\frac{\pi yj}{\pi jn^2} + \sum_{j=[4/a\pi]}^{[n/\pi]}\frac{\pi yj^2}{\pi jn^2\exp(cj)} + \sum_{j=[n/\pi]}^{n-1}\frac{\pi yj\exp(-C_2n)}{\pi jn}}\right)\\
& \leq & C\left(\frac{y}{n^2} + \frac{y}{n^2} \displaystyle{\sum_{j=[4/a\pi]}^{[n/\pi]}(\frac{j}{\exp(cj)})} + 4y\exp(-C_2n)\right) \leq C\frac{y}{n^2}\\
\end{eqnarray*}
\end{proof}

\begin{lem}\label{f(n,y)}If $f(x,y)$ is the solution of the Dirichlet problem in $\mathcal{R}(n)$ with $\phi(y) \equiv 1$, then there exists a constant $K>0$ such that for all $n$ and all $y \in \{1,..,n-1\}$, 
$$f(n,y) \leq K\frac{y}{n}.$$
\end{lem}

\begin{proof}

This is a particular case of Lemma \ref{Infrectangle}. From the proof of Lemma \ref{f(1,y)}, we know that
$$b_j \leq \frac{C}{j} \text{ and } \sin(\frac{\pi yj}{n}) \leq \frac{\pi yj}{n}.$$
We also know from Lemma \ref{a_jlemma} that for $1\leq j\leq n-1, a_j \geq \frac{j}{2n}$, so that
$$\exp(-a_jn) \leq \exp(-j/2).$$
%For $[n/\pi] \leq j \leq n-1, a_j \geq 1/2$, so that we also have
%$$\exp(-a_jn) \leq \exp(-n/2) \leq \exp(-j/2).$$
Plugging all this into (\ref{infrectangle}) gives
$$f(n,y) \leq K\displaystyle{\sum_{j=1}^{n-1}\frac{1}{j}\exp(-j/2)\frac{yj}{n}} \leq K\frac{y}{n}.$$
\end{proof}

\bibliographystyle{plain}
\bibliography{Estimates}

\end{document}